\documentclass[review,hidelinks]{siamart171218}

\usepackage[utf8]{inputenc}

\usepackage{amsfonts}
\usepackage{dsfont}  %
\usepackage{mathrsfs} 
\usepackage{graphicx}
\usepackage{comment}
\usepackage{graphicx}      

\usepackage{amsmath,amssymb}
 \usepackage{latexsym}
 \usepackage{multirow}

 \usepackage{bbm}

\usepackage{array}
\usepackage{makecell}
\usepackage{multirow}
\usepackage{xcolor}
\usepackage[T1]{fontenc}
\usepackage{eulervm}

\usepackage{hyperref} 
\usepackage{tikz}

\makeatletter
\DeclareRobustCommand{\sqcdot}{\mathbin{\mathpalette\morphic@sqcdot\relax}}
\newcommand{\morphic@sqcdot}[2]{%
  \sbox\z@{$\m@th#1\centerdot$}%
  \ht\z@=.33333\ht\z@
  \vcenter{\box\z@}%
}
\makeatother

\definecolor{purple}{RGB}{69,22,170}

\definecolor{verde}{RGB}{83,134,31}

\newcommand{\trp}{\intercal}

 \DeclareMathOperator{\Diag}{Diag}
 \DeclareMathOperator{\diag}{diag}
 \DeclareMathOperator{\tr}{tr}
 \DeclareMathOperator{\sign}{sign}
 \DeclareMathOperator{\cov}{Cov}
\DeclareMathOperator{\Ric}{Ric}

\ifpdf
  \DeclareGraphicsExtensions{.eps,.pdf,.png,.jpg}
\else
  \DeclareGraphicsExtensions{.eps}
\fi

\newsiamremark{assumptions}{Assumptions}

\newsiamremark{remark}{Remark}

\begin{nolinenumbers}
\ifpdf
\hypersetup{
  pdftitle={The $\mathcal{H}_2$-optimal Control Problem of CSVIU Systems: Discounted, Counter-discounted and Long-Run Solutions~Part II: Optimal Control},
  pdfauthor={Jo\~ao B. R. do Val and Daniel S. Campos}
}
\fi

\title{The $\boldsymbol{\mathcal{H}_2}$-optimal Control Problem of CSVIU Systems: Discounted, Counter-discounted and Long-Run Solutions\\Part II: Optimal Control%
\thanks{Submitted to the SIAM J Control  Optim editors on \today
\funding{Research supported in part by the National Council for Scientific and Technological Development (CNPq), grant n. 303352/2018-3, by FAPESP under grant n. 2016/08645-9, and by the Coordenação de Aperfeiçoamento de Pessoal de Nivel Superior - Brasil (CAPES) - Finance Code 001.}}}

\author{ Jo\~ao B. R. do Val\thanks{
UNICAMP, School of Electrical and Computer Engineering, Av. Albert Einstein 400, Campinas, SP, Brazil (\email{jbosco@fee.unicamp.br}, \email{d211498@dac.unicamp.br})}
\and
Daniel S. Campos\footnotemark[2]
}

\headers{The $H_2$-optimal Control of CSVIU Systems -- Part II}{Jo\~ao B. R. do Val and Daniel S. Campos}
\end{nolinenumbers}
\begin{document}
\bibliographystyle{siamplain}

\maketitle

\begin{nolinenumbers}
\begin{abstract}  
The paper deals with stochastic control problems associated with $H_2$ performance indices such as energy or power norms or energy measurements when norms are not defined. They apply to a class of systems for which a stochastic process conveys the underlying uncertainties, known as CSVIU (Control and State Variation Increase Uncertainty). These indices allow various emphases from focusing on the transient behavior with the discounted norm to stricter conditions on stability, steady-state mean-square error and convergence rate, using the optimal overtaking criterion---the long-run average power control stands as a midpoint in this respect.
 A critical advance regards the explicit form of the optimal control law, expressed in two equivalent forms.  One takes a perturbed affine Riccati-like form of feedback solution; the other comes from a generalized normal equation that arises from the nondifferentiable local optimal problem. They are equivalent, but the latter allows for a search method to attain the optimal law. 
 A detectability notion and a Riccati solution grant stochastic stability from the behavior of the norms. The energy overtaking criterion requires a further constraint on a matrix spectral radius.
With these findings, the paper revisits the emerging of the inaction solution, a prominent feature of CSVIU models to deal with the uncertainty inherent to poorly known models. Besides, it provides the optimal solution and the tools to pursue it.

\end{abstract}



\section{Introduction}

An exemplary mathematical model of a process of interest forms the ideal basis for forecasting future behavior, evaluating critical phenomena or measurements, or designing a control system.  Nevertheless, models are mostly nothing more than feeble imitations of reality, even though a good model can carry many of the essential features of the real world, \cite{bib:1,MKac}. As pointed out by many authors, this idea translates, in part, that to have more representative models for existing systems, one ought to characterize the uncertainties adequately,  c.f.  \cite{bib:11,bib:2}.  The need to deal with poor models is one of the main drives in the control field during the last three decades, which has brought a flourishing of new ideas that permeates many of the present date research, e.g., see \cite{FrancisKhar,bib:15,DraganMorozanStoica}.

Along these lines, the authors claim that the CSVIU model is a stochastic-based alternative to the uncertainty representation inherent to poor modeling of dynamic systems, cf.  \cite{bib:4,Fernandes2020,bib:3,CSVIU:fishery}. The CSVIU concept avoids the usual worst-case analysis of robust control by introducing a particular stochastic perturbation. A feature of the CSVIU method is a perturbing stochastic process that produces more significant drifts as the system state and control deviate from a better-known operation point, providing an interesting ground for accounting modeling errors.  The approach builds on the idea that when models are frails, state and control variations increase uncertainty, and cautious controls should prevail. The idea is to inbuild mathematically the notion that conservativeness should take place in the face of uncertainties. Solving the underlying optimal stochastic control problem, a region on the state space arises, in which an inaction control is optimal.  The zero-control (or zero-variation control) is the optimal feedback control therein, a behavior not encountered in the robust worst-case analysis, as far as the authors are aware. For more context and motivations, see the companion paper \cite{CSVIU:norm}.

This paper focuses on the $H_2$-norm control or the optimal overtaking control inspired on the norm measurement. The companion paper presents a study on the $H_2$-norm and associated energy and power measurements used here. 

We investigate the controlled version of the dynamical system studied in \cite{CSVIU:norm} by introducing a control variable $u(\cdotp) \in \mathds{R}^m$, as in \cite{bib:3}. 
Let us consider the following discrete-time process defined in a complete filtered probability space, $(\Omega, \mathcal{F}, P, \{ \mathcal{F}_k \}_{k \geq 0})$, the controlled CSVIU model, 
\begin{equation}
\Theta_{\textrm{ctr}} := 
\begin{cases}
\begin{aligned}[b]
        &x(k+1)  = A x(k) + ( \sigma_x + \overline{\sigma}_x \diag(|x(k)|)) \varepsilon^x (k) \\
        & \qquad\qquad + B u(k) + ( \sigma_u + \overline{\sigma}_u \diag(|u(k)|)) \varepsilon^u (k)
        + \sigma \omega (k),
    \end{aligned}\\
y(k) = Cx(k) + D u(k), \quad x(0) = x, 
\end{cases}
\label{CSVIUsystem_control}
 \end{equation}
 where $x\in\mathds{R}^n$, $A, \sigma_x, \overline{\sigma}_x \in \mathds{R}^{n \times n}$, $\sigma\in\mathds{R}^{n\times r}$ and $C\in \mathds{R}^{p \times n}$. In addition the control elements are  $B, \sigma_u, \overline{\sigma}_u \in \mathds{R}^{n \times m}$,  and $D\in \mathds{R}^{p \times m}$. The noise $\{\omega (k)\}_{k\ge0}$ is a $r$-dimensional persistent disturbance noise sequence, the ``nature noise'', and the ``system imperfection noise'' $\{ \varepsilon ^x(k),  \varepsilon ^u(k)\}_{k \geq 0}$ are $n$-dimensional and $m$-dimensional sequences, respectively. These are i.i.d. sequences with zero mean and their joint covariance forms an identity matrix. 
 
 The filtration $\{ \mathcal{F}_k \}_{k \geq 0}$ is composed by the sub-$\sigma$-algebras $\mathcal{F}_k \subset \mathcal{F}$ generated by the random variables $\omega(0),\varepsilon^x(0),\varepsilon^u(0), \ldots, \omega(k),\varepsilon^x(k),\varepsilon^u(k)$. Given a $n$-valued vector $x, \diag(|x|)$ is the $n$-dimensional diagonal matrix formed by setting $|x| = [\begin{smallmatrix} |x_1|& |x_2|&\cdots & |x_n|\end{smallmatrix}]^\trp$ as its diagonal, where $| \cdot |$ is the absolute value. The processes $\{x (k)\}_{k\ge0}$ is the $n$-dimensional state and $\{y (k)\}_{k\ge0}$ is a $m$-dimensional output of interest with $m\le n$.
All admissible control policies $\{u(k)\}_{k\ge0}$ belong to the set of $\{ \mathcal{F}_k \}$-adapted $m$-dimensional processes in feedback form, $\mathcal{U}:=\{u: u(t) = u_t (x_t ),{0\le t \le \kappa}\} $ with $\kappa$ finite or infinite. 

The difference equation in $\Theta_{\textrm{ctr}}$ is the controlled dynamics, whereas the other expresses the output variable.
The state variable and control inputs represent either the original variables or variations near an equilibrium point. Their absolute values (or absolute values variations) appear in the perturbation terms for the linear part of the model, representing measures of uncertainties associated with their departures from an equilibrium point, cf. \cite{bib:3}, or \cite{CSVIU:norm}.

For a control $u\in\mathcal{U}$, let us consider the following $\ell_2 (\Omega, \mathcal{F}, P)$ mean energy measurements for some $\alpha>0$,
\[
\mathcal{E}^{\kappa,\alpha}_2(u(\cdotp),y(\cdotp)) := E_{x} \Bigl[  \sum_{k=0}^{\kappa} \alpha^k \| y(k) \|^2 \Bigr], \quad\kappa\ge0
\]
in which $E_{x} \big[ \cdot \big]$ is the short for the expectation $E \big[ \cdot | \mathcal{F}_0 \big]$ and $y(\cdotp)$ is the corresponding output of $\Theta_{\textrm{ctr}}$ associated to the pair $(x(\cdotp),u(\cdotp))$. The optimal $H_2$-$\alpha$-energy-norm  is defined as
\begin{equation}
    \mathcal{E}^\alpha_2:=\inf_{u\in\mathcal{U}}\lim_{\kappa\to\infty}\mathcal{E}^{\kappa,\alpha}_2(u(\cdotp),y(\cdotp)), \quad y(0)=0 
    \label{energymeasure1}
\end{equation}    
 whenever finite. The pair  $(x^*(\cdotp),u^*(\cdotp))$ that attains the infimum of the energy, $\mathcal{E}^\alpha_2$, is called optimal. If $\alpha <1$ it is a discounted measurement, which is an adequate form to deal with the noise persistent excitation in system $\Theta_{\textrm{ctr}}$. When $\alpha\ge1$, we called it a counter-discounted measurement and \eqref{energymeasure1} is possibly unbounded and the norm is not defined. Nevertheless, energy measurements can still be compared by means of an overtaking criteria, cf. \cite{CarlsonHaurieLeizarowitz}, \cite{HuangYongZhou}. 
Optimality is expressed in the \emph{optimal overtaking  sense} as
\begin{equation}
   \mathcal{E}^{\alpha, \kappa}_2 (u^*(\cdotp),y^*(\cdotp))\le \mathcal{E}^{\alpha, \kappa}_2 (u(\cdotp),y(\cdotp))+\epsilon, \quad \forall u\in\mathcal{U}  
\end{equation}
for all $\kappa\ge\kappa_0$ sufficiently large. The output $y^*(\cdotp)$ corresponds  to the optimal $\Theta_{\textrm{ctr}}$-pair $(x^*(\cdotp),u^*(\cdotp))$.
 
 A third performance index is a power measurement based on the C\`esaro mean form,
\begin{equation*}
    \mathcal{P}_2(u(\cdotp),y(\cdotp)) := \lim_{\kappa \to \infty} \frac{1}{\kappa} E_{x} \Bigl[  \sum_{k=0}^{\kappa-1} \| y(k) \|^2 \Bigr]
\end{equation*}
provided that the limit exists. The optimal $H_2$-power-norm is defined as
\begin{equation}
    \mathcal{P}_2:=\inf_{u\in\mathcal{U}}\mathcal{P}_2(u(\cdotp),y(\cdotp)) 
    \label{energymeasure2}
\end{equation}

To accompany with these criteria, we introduce closely associated stochastic stability notions in \cite{CSVIU:norm}.
 Consider the  following notions, where $\mathcal{E}^\alpha_2(x(\cdotp)), \mathcal{E}^{\kappa,\alpha}_2(x(\cdotp))$ and $\mathcal{P}_2(x(\cdotp))$ are short for $\mathcal{E}^\alpha_2(0,x(\cdotp))$, $\mathcal{E}^{\kappa,\alpha}_2(0,x(\cdotp))$ and $\mathcal{P}_2(0,x(\cdotp))$, respectively.
 
\begin{definition}  \label{def:stability_definitions}
    System $\Theta_{\textrm{ctr}}$ is \emph{$\alpha$-stochastic stable}, 
    \begin{enumerate}
        \item[i)]  If $0\le\alpha<1$ and $\mathcal{E}^\alpha_2(x(\cdotp))=\lim_{\kappa\to\infty} \mathcal{E}^{\alpha,\kappa}_2(x(\cdotp))< \infty$,  $\forall x(0) = x_0 \in \mathds{R}^n$, 
        
        \item[ii)] If $\alpha\ge1$ and there exist $c_0,c_1$,  and $\xi_\kappa\in\mathds{R}^n$ with $\|\xi_\kappa\|\le c_2$   
        such that,  
        \[
        \mathcal{E}^{\kappa,\alpha}_2(x(\cdotp)) \le c_0\|x_0-\xi_\kappa\|^2+ \kappa c_1\alpha^\kappa,  \;\forall x(0) = x_0 \in \mathds{R}^n\text{ and $\kappa\ge0$,}
        \]
        
        \end{enumerate}
        
     System $\Theta$ is \emph{stochastic stable} if $\mathcal{P}_2(x(\cdotp))\le \bar{c} < \infty$, $\forall x(0) = x_0 \in \mathds{R}^n$.

\smallskip        
   The corresponding notion of stochastic stabilizability applies to system $\Theta_{\textrm{ctr}}$ whenever there exists some $u\in\mathcal{U}$ able to render one of the measurements, $\mathcal{E}^\alpha_2(u(\cdotp),x(\cdotp))$,  $\mathcal{E}^{\kappa,\alpha}_2(u(\cdotp),x(\cdotp))$ and $\mathcal{P}_2(u(\cdotp),x(\cdotp))$,  to satisfy the respective stability criterion.

\end{definition}

The link between the stability notions in Definition \ref{def:stability_definitions} and those output $H_2$ energy and power norms, or energy measurements criteria, is developed in \cite{CSVIU:norm} by introducing proper notions of stochastic detectability, see Definition 2.4 and Theorem 2.5 therein (here quoted as Definition \ref{def:detability.1} and Proposition \ref{prop:detectability}, respectively).

The discounted LQ control problem for a CSVIU model in discrete-time appears in \cite{bib:3}. To set up a link between the finiteness of the LQ problem with $\alpha$-stability, the authors impose the positivity condition $C^\trp C\succ0$. Here, with the results in \cite{CSVIU:norm} concerning $H_2$-norms and stability, stochastic detectability of $\Theta_{\textrm{ctr}}$ weakens that requirement. Moreover, when $\alpha\ge1$ we count on other conditions for $\alpha$-stability in \cite{CSVIU:norm}, associated with limited growth per stage of energy measurements as in Definition \ref{def:stability_definitions} (ii). With that, we are able to solve the control problem in the optimal overtaking sense, with assured stochastic stability and thus imposing severe requirements on steady-state behavior. As one sets $\alpha>1$, we can assess exponentially fast state and output convergence to a value in the mean-square sense, see \cite[lem.\,3.3,\;rem.\,3.4]{CSVIU:norm}.

A third novelty resides in the general form of optimal control. It can take two representations: one that involves a perturbed form of an affine Riccati solution, the other coined in terms of generalized gradients. 

We explore the two forms---the first stands as an optimal feedback global representation that arises from the local optimization problem, which is new to CSVIU models. With that, we can prove that the optimal feedback control is stabilizing for system $\Theta_{\textrm{ctr}}$ in the stochastic senses of Definition \ref{def:stability_definitions}. Equally importantly, the analysis of the second form provides the means to find the optimal solution, employing a nondifferentiable optimization method suited to the underlying problem with assured convergence.

\paragraph{The Inaction Region}
Because of the absolute values involved in the dynamical system model $\Theta_{\textrm{ctr}}$, an interesting feature of the optimal control for CSVIU is that it gives rise to a state space region in which the optimal control action is to sustain the action while the state variable is inside it. For each control component $u_i^*(x)$, $i \in \{ 1,2, \ldots , m \}$ of the optimal control vector $u^*(x) \in \mathds{R}^m$, we identify three partitions of state space $\mathds{R}^n$:
\begin{equation*}
\mathcal{R}_i^+ := \{ x \in \mathds{R}^n | u_i^*(x) > 0 \},\quad
    \mathcal{R}_i^- := \{ x \in \mathds{R}^n | u_i^*(x) < 0\},\quad
    \mathcal{R}_i^0 := \{ x \in \mathds{R}^n | u_i^*(x) = 0\}.
\end{equation*}

Each set $\mathcal{R}_i ^0$ may be a section of an $n$-dimensional region, and as a result, each $i$-component of the control solution subdivides the state space into disjoint regions.  The region $\mathcal{R}^0 = \bigcap_{i=1}^m \mathcal{R}_i^0$, which may be empty, is particularly called \emph{global inaction region} as all control inputs $ u_i$ are null.

This halt on the control action creates the notion of a cautionary inaction region and cautious controls. Accounting for the fact that the existing system is not precisely known, the CSVIU model originates optimal controls carrying the notion that in some situations, the best alternative is not to interfere in the system's dynamical evolution. The state vector then drifts uncontrolled while inside an optimization-designed neighborhood of the equilibrium point. This behavior is long seen in the finance and economic literature but, to the best of our knowledge, it does not appear to any other robust control approach in the system field, cf. \cite{bib:4}. 

On another extreme, consider regions of state space far enough from $\mathcal{R}_i ^0$ such that $u_i ^*(x) \neq 0$ for each $i = 1,\cdots,m$. The optimal control signal vector is $s_u := [\begin{smallmatrix}s_{u,1}& \cdots& s_{u,m}\end{smallmatrix}]^\trp$ with each $s_{u,i}=\pm1$ whenever 
$u_i^*(x) \gtrless 0 $, and \emph{homogeneous signals regions} are defined by $\mathcal{R}^{s_u} = \cap_{i =1}^m \mathcal{R}_i ^{s_{u,i}}$ for each $s_u$. Those are  state space regions 
in which the signal vector $s_u$ of $u^*(x)$ remains unaltered. 
Since changes of the state vector signals also impact the solution, one keeps track of the open orthants sets,  $\mathbb{O}_j, j= 1,\ldots 2^n $ of $\mathds{R}^n$. In this scenario the idea of \emph{asymptotic regions} emerges,  namely, $\mathds{R}^n$ regions for which state and the optimal control vectors are nonzero with given signals, i.e.,   $x \in \mathcal{R}^{s_u}\cap \mathbb{O}_j$, and such that $\| x \| \to \infty$ following directions leading far away from any of the signal switching boundaries.

That is the framework laid so far for CSVIU in discrete-time \cite{bib:3} and for Brownian motion driven continuous-time  \cite{bib:4}. Here we report considerable advances in the form of the optimal control global law for $H_2$-norms and optimal overtaking problems in discrete-time. The qualitative behavior just described is made redundant by the present analysis. 

Stochastic detectability and the solution of a perturbed algebraic Riccati equation play an crucial role in devising optimal stabilizing feedback laws for system $\Theta_{\textrm{ctr}}$. The importance of such a Riccati equation was previously noted, connected with the aforementioned asymptotic regions, hence restricted to a local representation. Here, we reveal the general optimal solution and put the Riccati equation role into the proper perspective. Besides, we refine the notion of the inaction regions and under the light of the results, develop a nondifferentiable optimization method to deal with the local optimization problem. The solution sheds some light on the shaping of the inaction region. The method has assured convergence and can solve numerically the optimal control problem engendered here.

 The paper is organized as follows. Section \ref{sec:analysisH2norm_controlled} quotes the results of \cite{CSVIU:norm} concerning energy measurements and norms, stochastic stability and detectability, and closes with two auxiliary lemmas. 
 Section \ref{sec:optimal.solution} refers to the solution of the CSVIU $H_2$ energy (discounted) and power (C\'{e}saro average) optimal controls, together with the optimal overtaking (counter-discounted) criterion and their stabilizing property 
 Section \ref{se:inaction.representation} explores the optimal solution from the nondifferentiable form. It completes the solution construct, the inaction region and helps devise a convergent method to obtain the solution. Finally, Section \ref{sec:conclusion} gathers some concluding remarks.

\section{Notation and Preliminary Results}
\label{sec:analysisH2norm_controlled}

For a matrix $U\in\mathds{R}^{n\times n}$, $\Diag(U) \in \mathds{R}^{n \times n}$ is the diagonal matrix formed by the main diagonal of $U$, $U_d\in\mathds{R}^n$ is the vector in the main diagonal of $U$ and $\tr(U)$ denotes the trace operator.  $\mathds{S}^{n+}$ stands for the real vector space of symmetric matrices of size $n$ that are positive semidefinite. For  $U\in\mathds{S}^{n+}$, $U \succ 0$ ($U\succeq0$) indicates that $U$ is a positive (semi-) definite matrix, also for $U,V\in\mathds{S}^{n+}$, $V\succeq U \Leftrightarrow V-U\succeq 0$. $\|U\|$ indicates any matrix norm and for a square matrix $\lambda^+(U) (\lambda^-(U))$ denotes the largest (minimum) eigenvalue of $U$ and $r_\sigma(U)$ its spectral radius.
 
For a vector $u \in \mathds{R}^n$, $|u|$ indicates the vector $\begin{bmatrix}
|u_1|&|u_2|&\cdots&|u_n|
\end{bmatrix}^\trp$ and $\diag(u)$ stands for the diagonal matrix made up by vector $u$.
For two vectors $u,v \in \mathds{R}^n$, $\langle u , v \rangle$ denotes the inner product,  $u \sqcdot v$ denotes the Hadamard product, and the square (semi-)norms $\|x\|^2_U$ stands for  $\langle x, U x\rangle$, in which $U \in \mathds{S}^{n+}$.

In connection with the data of system $\Theta_{\textrm{ctr}}$, let us define the linear operators $\mathcal{Z}_x: \mathds{S}^{n+} \to\mathds{S}^{n+}$, $\mathcal{W}_x: \mathds{S}^{n+} \to \mathds{S}^{n}$, given by:
\begin{subequations}\label{eq:state-operators} 
\begin{align}\allowdisplaybreaks
    \mathcal{Z}_x(U) &= \Diag ( \overline{\sigma}_{x}^{\trp} U \overline{\sigma}_x ),
    \label{linear_operators_definition}
\\
\mathcal{W}_x(U)  &= \Diag(\overline{\sigma}_x^\trp U \sigma_x + \sigma_x^\trp U \overline{\sigma}_x ),
     \label{linear_operators_definition.2}
\intertext{together with $\mathcal{L}^\alpha:\mathds{S}^{n+} \to \mathds{S}^{n+}$, for some $\alpha \geq 0$,}
    \mathcal{L}^\alpha(U) &= \alpha(A^\trp U A + \mathcal{Z}_x(U)), 
    \label{special_linear_operators_definition}
\intertext{and sometimes we denote $\mathcal{L}^1(U) = \mathcal{L}(U)$.
In addition to operators in \eqref{eq:state-operators}, let us adopt counterparts $\mathcal{Z}_u: \mathds{S}^{n+} \to \mathds{S}^{n+}$, $\mathcal{W}_u: \mathds{S}^{n+} \to \mathds{S}^{n}$, and $\varpi_1 : \mathds{S}^{n+} \to \mathds{R}$ such that for any $U \in \mathds{S}^{n+}$,}
    \mathcal{Z}_u(U) & = \Diag(\overline{\sigma}_u^\trp U \overline{\sigma}_u),\\
    \mathcal{W}_u(U) & = \Diag(\overline{\sigma}_u^\trp U \sigma_u + \sigma_u^\trp U \overline{\sigma}_u),\\
    \varpi_1(U) & = \tr \{ U(\sigma \sigma^\trp + \sigma_x \sigma_x^\trp + \sigma_u \sigma_u^\trp ) \},
\label{linear_operators_definition2}
\intertext{and $\Sigma: \mathds{R}^{n \times n} \to \mathds{R}^{m \times n}$, $\Lambda: \mathds{R}^{n \times n} \to \mathds{R}^{m \times m}$ with}
    \Sigma(U) &=  B^\trp U A + \alpha^{-1}D^\trp C,
    \label{rewriteitoninfinitynorm1}
\\
    \Lambda(U) &= B^\trp U B + \mathcal{Z}_u(U) + \alpha^{-1} D^\trp D.
    \label{rewriteitoninfinitynorm2}
\intertext{When $\Lambda$ is invertible, consider also $\Ric^\alpha : \mathds{R}^{n \times n} \to \mathds{R}^{n \times n}$ such that}
\Ric^\alpha(U) &=\mathcal{L}^\alpha(U)  - \alpha\Sigma(U)^\trp \Lambda(U)^{-1}\Sigma(U) + C^\trp C
\label{eq:riccati.oper}
\end{align}
\end{subequations}

With the exception of $\mathcal{W}_x$ and $\mathcal{W}_u$ these are all linear-positive operators, i.e., $U \succeq 0$ implies $\Pi(U) \succeq 0$, where $\Pi$ stands for any of the operators in \eqref{eq:state-operators}. 
A linear-positive operator is also monotone, namely, if $U \succeq V$ for $U,V \in \mathds{S}^{n+}$ then $\Pi(U)\succeq\Pi(V)$. Here,  $U\succeq V \Leftrightarrow U - V\succeq 0$.
When convenient, we refer by $\mathcal{W}_{u_d}$ (or $\mathcal{W}_{u_d}(U)$) to the $\mathds{R}^m$-vector in the diagonal matrix  $\mathcal{W}_{u}$. Similarly,
with $\mathcal{W}_{x_d}$ extracted from $\mathcal{W}_{x}$.

Let us consider the set $\{ -1, 0, +1\}$ and create the collection $\mathscr{S}= \{ -1, 0, +1\}^n$ of distinct vectors $s_i \in\mathscr{S}, i= 1, \ldots , 3^n $ under some arbitrary order.
We introduce the signal vector function $\mathcal{S}: \mathds{R}^n \to \mathscr{S}$ such that for each $x\in\mathds{R}^n$,
\begin{equation}
   \mathcal{S}(x) = \begin{bmatrix}\sign(x_1) & \cdots & \sign(x_n)\end{bmatrix}^\trp
   \label{eq:signalfunction}
\end{equation}
with the convention $\sign(0)=0$.

\paragraph{Notions of $\alpha$-stability and detectability}
For ease of reference in what follows, we quote essential results of \cite{CSVIU:norm}. We refer to system $\Theta$, the uncontrolled version of $\Theta_{\textrm{ctr}}$, namely, $u(\cdotp)\equiv0$, and we drop the control dependence in the notation. To recast the same ideas to the controlled setting, readily recall the notions of stochastic stabilizability introduced in Definition \ref{def:stability_definitions}.

 The starting point is the study of positive operators in ordered Banach spaces applied to linear perturbed Lyapunov equations that arise in many stochastic problems, e.g., \cite{freiling,Hasanov}. We deal with the linear-positive operators listed in \eqref{eq:state-operators}.

\begin{proposition}[prop.\;2.2 \cite{CSVIU:norm}]\label{prop:freiling.stab1}
 $\alpha$-stability of $\Theta$ for some $0\le\alpha<1$ is equivalent to require that
\begin{enumerate}
    \item[i)] $\mathcal{L}^\alpha$ is an inverse-positive operator,
    \item[ii)] $\mathcal{L}^\alpha$ is $d$-stable,
    \item[iii)] There exists $U\succ0$ such that $(I-\mathcal{L}^\alpha)(U)\succ0$,
    \item[iv)] $\alpha^{1/2}A$ is $d$-stable relative to $\alpha \mathcal{Z}_x$,
    \item[v)] All eigenvalues of $\sqrt{\alpha}A$ lie in the open unit disk and $r_\sigma((I-\alpha\mathds{A})^{-1}\mathcal{Z}_x)<\alpha^{-1}$, where $\mathds{A}(U) := A^\trp U A$ for $U\in\mathds{R}^{n\times n}$.
\end{enumerate}

If $\alpha\ge1$ and all eigenvalues of ${\alpha}A$ lay in the open unit disk then (i)--(v) are equivalent to $\alpha$-stochastic stability of $\Theta$.

If all eigenvalues of $A$ lay in the open unit disk then  (i)--(v) with $\alpha=1$  are equivalent to stochastic stability of $\Theta$.
\end{proposition}

\begin{definition}[$\alpha$-detectability]\label{def:detability.1}
 System $\Theta$ is $(C,\alpha^{1/2}A,\alpha\mathcal{Z}_x)$ $\alpha$-detectable if there is a matrix $H \in \mathds{R}^{n \times p}$ such that $\alpha^{1/2}(A+HC)$ is $d$-stable relative to $\alpha\mathcal{Z}_x$. 
Most of the times, this notion is referred as $(C,\mathcal{L}^\alpha)$-detectable for some $\alpha\ge0$, or as $(C,\mathcal{L})$-detectable when $\alpha=1$.
\end{definition}

 The notion above and Proposition \ref{prop:freiling.stab1} lead to the result.

\begin{proposition}[thm.\;2.5\,\cite{CSVIU:norm}]\label{prop:detectability}
Suppose that for some $\alpha\ge0$, $(I-\mathcal{L}^\alpha)(U) =C^\trp C$ has a solution $L\succeq0$. If $\Theta$ is $(C,\mathcal{L}^\alpha)$-detectable then if either,
\begin{enumerate}
    \item[i)] $\alpha<1$, or,
    \item[ii)] $\alpha\ge1$ and all eigenvalues of $\alpha A$ lies in the open unit disk,
\end{enumerate}
then  $\Theta$ is $\alpha$-stochastically stable and conversely. 
 \begin{enumerate}
     \item[iii)] When (ii) holds for $\alpha=1$, $\Theta$ is stochastically stable and conversely. 
 \end{enumerate}
\end{proposition}

\paragraph{Controlled System}
Returning to system $\Theta_{\textrm{ctr}}$, we adopt the compact notation,
\begin{equation}\label{eqs:compact.theta}
\begin{aligned}[c]
       x_{k+1} &= A x_k + Bu_k + \sigma_x(x_k)\omega_{0,k}+ \sigma_u(u_k)\epsilon^u_{k}
     \\
      y_k &= Cx_k+Du_k, \quad k\ge0
\end{aligned}
 \end{equation}
 such that 
 $   \sigma_x(x)  = \begin{bmatrix}
     \sigma & \sigma_x + \bar{\sigma}_x \diag(|x|)
    \end{bmatrix}
$, $\omega_{0,k} = [ \omega(k) ~~ \varepsilon^x(k)  ]^\trp$  and $\sigma_u(u) = \sigma_u + \bar{\sigma}_u \diag(|u|)$, $\epsilon^u_{k} = \epsilon^u(k)$. Recall that  $\cov([\omega_{0,k}^\trp \;(\epsilon^u_k)^\trp])=I_{r+n+m}$. The associated processes and time stage sequences are  expressed using the stage $k$ as a subindex.

Let us consider some sequences $\{P_k\}, P_k\in\mathds{S}^{n+}, \{r_k\}, r_k\in\mathds{R}^n$ and $\{g_k\}, g_k\ge0$, $k=0,1,\ldots$ and introduce functions $V^u:\mathbb{N}\times \mathds{R}^n\to \mathds{R}$ made up as, 
\begin{equation}\label{eq:Vu}
    V^u(k, x) := \alpha^k (x^\trp P_k x + \langle r_k , |x| \rangle + g_k),\quad x\in\mathds{R}^n.
\end{equation}
The proof of the next lemma appears in Appendix \ref{app:lemm.variation.Vu}.
 \begin{lemma}\label{lemm:variation.Vu}
 For any pair $k\to(x_k,u_k)$ of state and admissible control of $\Theta_{\textrm{ctr}}$,
  \begin{multline}\label{eq:dif.one.step.expected}\allowdisplaybreaks
     E\left[\alpha^{-k}(V^u_{k+1}( x_{k+1}) - V^u_k (x_k))+\|y_k\|^2 |x_k\right]=
    \\
      x_k^\trp \big(\mathcal{L}^\alpha(P_{k+1})+ C^\trp C -  P_{k} \big) x_k + \alpha u_k^\trp \Lambda(P_{k+1}) u_k +
  \\
       + \langle \alpha A^\trp \mu_{k+1} + \alpha \mathcal{W}_{x_d}(P_{k+1})\mathcal{S}(x_k) - \mu_{k} ,x_{k} \rangle 
    \\
        +\alpha\langle  B^\trp \mu_{k+1} + \mathcal{W}_{u}(P_{k+1})\mathcal{S}(u_k) +2\Sigma(P_{k+1})x_k ,u_{k} \rangle  
         +\alpha g_{k+1} + \alpha \varpi_1(P_{k+1}) - g_{k}
\end{multline}
where $\mu_k:=E[r_k\sqcdot\mathcal{S}(x_k)|x_k]$ and $\mu_{k+1}:=E[r_{k+1}\sqcdot\mathcal{S}(x_{k+1})|x_k]$.
 \end{lemma}

For any $G\in\mathds{R}^{m\times n}$, let us denote $\mathcal{H}^\alpha_G:\mathds{S}^{n+} \to \mathds{S}^{n+}$  the operator,
\begin{equation}
\mathcal{H}^\alpha_G(U) := \mathcal{L}^\alpha(U)+ \alpha 
\begin{bmatrix}
I\\G
\end{bmatrix}^\trp
\begin{bmatrix}
\alpha^{-1} C^\trp C &\Sigma(U) 
\\
\Sigma(U)^\trp & \Lambda(U)
\end{bmatrix}
\begin{bmatrix}
I\\G
\end{bmatrix}    
 \end{equation}

Next, we set a parallel to the results on $\alpha$-stability in \cite{CSVIU:norm} for a system $\Theta_{\textrm{ctr}}$ that is $\alpha$-stabilizable by a linear feedback control. The proof appears in Appendix \ref{app:stabilization.G}.
\begin{lemma}\label{lemm:stabilization.G}
For $G\in\mathds{R}^{n\times m}$ define the sequences $P_k, \mu_k$ and $g_k$, for $k=0,\ldots,\kappa$, as 
  \begin{subequations}\label{eqs:EqDiffs.controled.G}
  \begin{gather}\label{eq:quadratic.G}
 \mathcal{H}_G^\alpha(P_{k+1})= P_k,  
\\\label{eq:linear.G}
  (A+B G)^\trp \mu_{k+1}  +  \mathcal{W}_{x}(P_{k+1})\mathcal{S}(x_k) + G^\trp\mathcal{W}_{u}(P_{k+1})\mathcal{S}(u_k) = \alpha^{-1} \mu_{k},
\\
\label{eq:scalar.G}
   g_{k+1} +  \varpi_{\textcolor{blue}{1}}(P_{k+1}) = \alpha^{-1} g_{k}
\end{gather}
  \end{subequations}
with $P_\kappa=0$, $\mu_\kappa=0$ and $g_\kappa=0$. Then for $\Theta_{\textrm{ctr}}$ with $k\to u_k=Gx_k$,
\begin{equation*}
   E\Bigl[\sum_{k=0}^{\kappa-1} \alpha^{k} \| y(k) \|^2\Bigr] = \|x_0\|^2_{P_0} + \langle \mu_0, x_0\rangle +g_0
\end{equation*}
Moreover, suppose that  $(C,\mathcal{L}^\alpha)$ is $\alpha$-detectable,
\begin{equation}\label{eq:controled.expanded}
 (I-\mathcal{H}_G^\alpha)(U) =0    
 \end{equation}
has a solution $L\succeq0$ and if $\alpha >1$ then $r_\sigma(A+BG)<\alpha^{-1}$. In this situation, $\Theta_{\textrm{ctr}}$ is $\alpha$-stabilizable and $k\to u_k=Gx_k$ $\alpha$-stabilizes $\Theta_{\textrm{ctr}}$.
 
\end{lemma}

\section{The Optimal Solution} \label{sec:optimal.solution}
The next theorem is in position to present the optimal solution. The first part of the proof develops the optimal law, the second shows the optimal norms and the optimal overtaking solutions. The third part proves that the optimal solution $\alpha$-stabilizes $\Theta_{\textrm{ctr}}$.

\begin{theorem}
 \label{theo:controlled.solution}
 Suppose that $D^\trp D\succ0$ and for some $\alpha> 0$ $\Theta_{\textrm{ctr}}$ is $(C,\mathcal{L}^\alpha)$-detectable, there is a solution $L\succeq0$ to $(I-\Ric^\alpha)(U) = 0$ and if $\alpha>1$,   $r_\sigma(A-B\Lambda(L)^{-1}\Sigma(L))<{1}/{\alpha}$ holds.
 For a  pair $k\to (x_k,u_k)$ of system $\Theta_{\textrm{ctr}}$, consider the difference equations for $k\to v_k$ and $k\to g_k, k=0,1,\ldots$, 
\begin{gather}\label{eq:reccurrence.mu}\allowdisplaybreaks
      (A-B\Lambda(L)^{-1}\Sigma(L))^\trp \; v_{k+1} +    
      \mathcal{W}_{x}(L) \mathcal{S}(x_k) - \Sigma(L)^\trp\Lambda(L)^{-1}\mathcal{W}_{u}(L)\mathcal{S}(u_k) \!=\! \frac{1}{\alpha}v_k,
\\
\label{eq:stage.residual}
\text{and~ }\alpha g_{k+1} + \alpha \varpi_{1}(L)  +  \rho_k = g_{k},    
\end{gather}
 with $\lim_{k\to\infty}v_k,g_k$ arbitrary, and where,
\begin{subequations}\label{eq:rho.min}
\begin{gather}
    \label{eq:stage.optimization}
     \rho_k:= \min_{u\in\mathds{R}^m} \rho_k(u),
\\
\label{eq:rho.L}
\begin{array}{lrl}
 \textrm{with}   &\rho_k(u)&\!\!\!:=\|u-u^0_k\|^2_{\Lambda(L)}-\dfrac{\alpha}{4}\left\|B^\trp \mu_{k+1} + \mathcal{W}_{u}(L) \mathcal{S}(u) \right\|^2_{\Lambda(L)^{-1}},\\[1.5ex]
    &u^0_k &\!\!\!:=   - \Lambda(L)^{-1}
\Bigl( \Sigma(L)x_k+ \dfrac{1}{2}\bigl(B^\trp {\mu}_{k+1} + \mathcal{W}_{u}(L) \mathcal{S}(u)\bigr)\Bigr).
\end{array}
\end{gather}
\end{subequations}
and $\mu_{k+1}=E[v_{k+1}|x_k]$.  Then, 
\begin{enumerate}
    \item[I.] For any output $k\to y_k$ induced by a pair $k\to (x_k,u_k)$,
\begin{equation} \label{eq:optimality.sense1}
  \|x_0\|^2_L + \langle E[\mu_{0}|x_0],x_0\rangle 
\le \sum_{k=0}^\infty\alpha^k E\bigl[\|y_k\|^2 -(\rho_k +\alpha \varpi_1(L)) \bigr]
\end{equation}
holds. Moreover, $k\to \rho_k$ is bounded, and if
\begin{equation}\label{eq:argument.optimization}
     \bar{u} _k =  \arg\min_{u\in\mathds{R}^m} \rho_k(u),\quad k\ge0
\end{equation}
then $k\to (\bar{x}_k,\bar{u}_k)$ is the optimal pair and the equality in \eqref{eq:optimality.sense1} is attained. 

\item[II.]For $\alpha\le1$, the norms are
\begin{equation}\label{eq:optomal.norms}
    \begin{cases}
     \mathcal{E}^\alpha_2 = \dfrac{\alpha}{1-\alpha}\varpi_1(L) +\sum_{k=0}^\infty \alpha^kE[\rho_k|x_0=0] , \quad\text{if $0\le\alpha<1$,}
    \\
     \mathcal{P}_2 =\varpi_1(L) + \lim_{k\to\infty} E[\rho_k].
    \end{cases}
\end{equation}
When $\alpha>1$, for each $\epsilon>0$, 
\begin{equation}\label{eq:overtaking.optimality}
     \mathcal{E}^{\alpha, \kappa}_2 (y(\cdotp))- \mathcal{E}^{\alpha, \kappa}_2 (\bar{y}(\cdotp))>  -\epsilon,
\end{equation}
holds for any $\kappa\ge\kappa_0$ sufficiently large, and an output $k\to y_k$ induced by any pair $k\to (x_k,u_k)$.

\item[III.] The feedback law $\{\bar{u}_k \}_{k=0,\ldots}$ in \eqref{eq:argument.optimization}  $\alpha$-stabilizes $\Theta_{\textrm{ctr}}$.
\end{enumerate}

\end{theorem}

The  proof of the theorem is divided into the corresponding assertions (I)--(II).

\begin{proof}[Proof of part (I)]
Let us consider functions $V^u_k$ in \eqref{eq:Vu}, and their expected variation \eqref{eq:dif.one.step.expected} in Lemma \ref{lemm:variation.Vu} for each $k$. Denote  $s_k:=\mathcal{S}(x_k), s_{k+1}:=\mathcal{S}(x_{k+1})$ and $s_u:=\mathcal{S}(u_k)$, and after some algebraic manipulations, we get
\begin{multline*}\allowdisplaybreaks
     E\left[\alpha^{-k}(V^u_{k+1}( x_{k+1}) - V^u_k (x_k))+\|y_k\|^2 |x_k\right]
    \\
      =x_k^\trp \big(\mathcal{L}^\alpha(P_{k+1}) + C^\trp C -  P_{k} \big) x_k 
       +\langle \alpha  A^\trp {\mu}_{k+1} + \alpha \mathcal{W}_{x}(P_{k+1}) s_k - \mu_k ,x_{k} \rangle 
      \\
      - \frac{\alpha}{4}\left\| B^\trp {\mu}_{k+1} + \mathcal{W}_{u}(P_{k+1}) s_u + 2\Sigma(L)x_k\right\|^2_{\Lambda(P_{k+1})^{-1}}
      \\
       + \|u_k-u^0_k\|^2_{\Lambda(P_{k+1})}+ \alpha g_{k+1} + \alpha \varpi_{\textcolor{blue}{1}}(P_{k+1}) - g_{k} 
\end{multline*}
where, under the assumptions, the inverse matrix exists, and we write
\begin{equation*}
u^0_k= 
\\
- \Lambda(P_{k+1})^{-1}
\Bigl( \Sigma(P_{k+1})x_k+ \frac{1}{2}\bigl(B^\trp {\mu}_{k+1} + \mathcal{W}_{u}(P_{k+1}) s_u\bigr)\Bigr)
\end{equation*}
Hence,
\begin{multline}\label{eq:final.adjust.control}\allowdisplaybreaks
           E\left[\alpha^{-k}(V^u_{k+1}( x_{k+1}) - V^u_k (x_k))+\|y_k\|^2 |x_k\right]=
      x_k^\trp \big(\Ric^\alpha(P_{k+1}) -  P_{k} \big) x_k 
      \\
      +\alpha\langle  (A-B\Lambda(P_{k+1})^{-1}\Sigma(P_{k+1}))^\trp {\mu}_{k+1},x_{k} \rangle
       \\
      +\langle\alpha\bigl( \mathcal{W}_{x}(P_{k+1}) s_k -
      \Sigma(P_{k+1})^\trp\Lambda(P_{k+1})^{-1}\mathcal{W}_{u}(P_{k+1})s_u  \bigr) -  \mu_k ,x_{k} \rangle 
       \\
       -\frac{\alpha}{4}\left\|B^\trp \mu_{k+1} + \mathcal{W}_{u}(P_{k+1}) s_u\right\|^2_{\Lambda(P_{k+1})^{-1}}
      \\
       + \|u_k-u^0_k\|^2_{\Lambda(P_{k+1})} + \alpha g_{k+1} + \alpha \varpi_{\textcolor{blue}{1}}(P_{k+1}) - g_{k} 
\end{multline}

Now, set $\varrho^*_k=\min_{u\in\mathds{R}^m} \varrho_k(u)$ where
\begin{subequations}\label{eqs:opt.choice}
\begin{align}\label{eq:rho.u}
&\varrho_k(u) :=  \|u-u^0_k\|^2_{\Lambda(P_{k+1})} 
     -\frac{\alpha}{4}\left\|B^\trp \mu_{k+1} + \mathcal{W}_{u}(P_{k+1}) \mathcal{S}(u) \right\|^2_{\Lambda(P_{k+1})^{-1}}, 
\\\label{eq:ricc.evolution}
&\Ric^\alpha(P_{k+1})=P_{k},
\\\label{eq:abs.evolution.controlled}
    &(A+BG_{k+1})^\trp \mu_{k+1} 
     +  \mathcal{W}_{x}(P_{k+1}) \mathcal{S}(x_k) + G_{k+1}^\trp\mathcal{W}_{u}(P_{k+1})\mathcal{S}(u)
     = \frac{1}{\alpha}\mu_k
\\
&\alpha g_{k+1} + \alpha \varpi_1(P_{k+1})  +  \varrho^*_k = g_{k} 
\end{align}
\end{subequations}
with $G_k= - \Lambda(P_{k})^{-1}\Sigma(P_{k})$. Then
$ E\big[ V^u_{k+1}( x_{k+1}) - V^u_k (x_k) +\alpha^k\|y_k\|^2|x_k\big]\ge0 $ holds true, and the equality is attained if one applies 
\begin{equation}\label{eq:optimal.law}
u^*_k =\arg\min_{u\in\mathds{R}^m} \varrho_k(u)
\end{equation}

Let us consider the feedback law $\{u^*_k\}_{ k=0,\ldots,\kappa}$ for $\kappa>0$, defined by \eqref{eq:optimal.law}. By setting $P_\kappa=0,  \mu_\kappa=0$ and $g_\kappa=0$, readily $V_\kappa^{u^*}(x) =0$. 
Since the process $\Theta_{\textrm{ctr}}$ is Markovian, and taking into account \eqref{eq:final.adjust.control}--\eqref{eqs:opt.choice},
\begin{multline}\label{eq:optimal.finite.cost.comp}
 V^{u^*}_0(x_0)  = 
   E[V^{u^*}_0(x_0) - V^{u*}_{\kappa}(x_\kappa) | x_0] 
   = 
   \\
   E\bigl[\sum_{k=0}^{\kappa-1}  E\left(  V^{u^*}_k (x_k) -V^{u^*}_{k+1}( x_{k+1}) |x_k\right)\big|x_0\bigr]=
   E\bigl[\sum_{k=0}^{\kappa-1}\alpha^k\bigl(\|y_k\|^2+ \rho^*_k - \rho_k(u_k)\bigr)|x_0\bigr]
    \\
    \le E\bigl[\sum_{k=0}^{\kappa-1}\alpha^k\|y_k\|^2|x_0\bigr]= \mathcal{E}_2^{\alpha,\kappa}(y(\cdot)),
\end{multline}
   is valid for any feedback law $\{u_k\}_{k=0,\ldots,\kappa}$. The equality is attained however, when $\{u^*_k\}_{k=0,\ldots,\kappa}$ as in  \eqref{eq:optimal.law} is applied. In other words, 
  \begin{equation}\label{eq:optimalcost.kappa}
  V_0^{u^*}(x_0)\le \mathcal{E}_2^{\alpha,\kappa}(y(\cdot))
  \end{equation}
and the equality is attained for  $\{u^*_k\}_{k=0,\ldots,\kappa}$.

Now, we explicitly indicate the horizon $\kappa$ as $P_k^{(\kappa)}$ or $\mu_k^{(\kappa)},\forall k<\kappa$ for the solutions of \eqref{eq:ricc.evolution} and \eqref{eq:abs.evolution.controlled} with $P_\kappa^{(\kappa)}=0$ and  $\mu_\kappa^{(\kappa)}=0$, respectively. From the assumptions we get that, $0\preceq P_k^{(\kappa)}\uparrow L$ in the semipositive definite sense, as $\kappa\to\infty$, where $L$ is the solution of $(I-\Ric^\alpha)(U)=0$. This is equivalent to say that $L\succeq0$ satisfies $ (I -\mathcal{H}^\alpha_G)(U)=0$ with $G= - \Lambda(L)^{-1}\Sigma(L)$, and since $(C,\mathcal{L}^\alpha)$ is $\alpha$-detectable it implies that  $\mathcal{A}:=A+BG$ is such that  $r_\sigma(\mathcal{A})<\alpha^{-1/2}$. If $\alpha >1$ an assumption requires that  $r_\sigma(\mathcal{A})<\alpha^{-1}$. 
As result, 
\begin{multline}\label{eq:vk.infty}
  \mu_k^{(\infty)} = \lim_{\kappa\to\infty}\mu^{(\kappa)}_k
  \\
  =\lim_{\kappa\to\infty} E\Bigl[\sum_{n=0}^\kappa\alpha^{n+1} (\mathcal{A}^\trp)^n\bigl(\mathcal{W}_{x}(P_{k+n}^{(\kappa)}) \mathcal{S}(x_{k+n})
+G^\trp\mathcal{W}_{u}(P_{k+n}^{(\kappa)})\mathcal{S}(u^*_k)\bigr)\Big| x_k\Bigl]
\\
=E\Bigl[\sum_{n=0}^\infty \alpha^{n+1} (\mathcal{A}^\trp)^n\bigl(\mathcal{W}_{x}(L) \mathcal{S}(x_{k+n})
+G^\trp\mathcal{W}_{u}(L)\mathcal{S}(u^*_k)\bigr)\Big| x_k\Bigl],\quad\forall k\ge0 
\end{multline}
and therefore, $  |\mu_k^{(\infty)}| \le \bar{\mu}$ where,
\begin{equation}\label{eq:bar.v.ctrl}
 \bar{\mu} = \alpha r_\sigma\left((I-\alpha \mathcal{A}^\trp)^{-1}\right)\Bigl(\bigl|\mathcal{W}_{x_d}(L)\bigr|
+\|G \| \bigl|\mathcal{W}_{u_d}(L)\bigr|\Bigr)
 \end{equation}
hence $k\to\mu_k^{(\infty)}$ is a bounded $n$-valued processes. Convergence of $\mu_k^{(\kappa)}$ and monotone convergence of each $P_k^{(\kappa)}$ lead to  convergence of $\varrho_{k}(u)$ in \eqref{eq:rho.u} to $\rho_{k}(u)$ in \eqref{eq:rho.L} (and thus, $\varrho^*_k \to\rho_k$ in \eqref{eq:stage.optimization}).
 Note that for $\rho_k$ in \eqref{eq:rho.min} that
 \begin{multline}
     \rho_k\ge  -\frac{\alpha}{4}\left\|B^\trp\mu_{k+1} + \mathcal{W}_{u}(L) \mathcal{S}(u) \right\|^2_{\Lambda(L)^{-1}}
     \\
     \ge -\frac{\alpha}{4}\left(\|B^\trp\|^2_{\Lambda(L)^{-1}}\|\bar{\mu}\|^2_{\Lambda(L)^{-1}} + \|\mathcal{W}_{u_d}(L)  \|^2_{\Lambda(L)^{-1}}\right)
 \end{multline}
In addition, if we set $u_k = G x_k$ as above, we get that
 \begin{multline}
     \rho_k \le \rho_k(Gx_k) \le \|Gx_k-u^0_k\|^2_{\Lambda(L)} 
     = \frac{1}{4}\left\|B^\trp \mu_{k+1} +  \mathcal{W}_{u}(L) \mathcal{S}(Gx_k) \right\|_{\Lambda(L)^{-1}}^2 
     \\
      \le  \frac{1}{4}\left(\|B^\trp\|^2_{\Lambda(L)^{-1}}\|\bar{\mu}\|^2_{\Lambda(L)^{-1}}   + \|\mathcal{W}_{u_d}(L)  \|_{\Lambda(L)^{-1}}^2 \right) 
 \end{multline}
 showing that $\rho_k$ is bounded.
\end{proof}

\begin{proof}[Proof of part (II)]
For the norm evaluation, consider first the case when $\alpha<1$. Set $P_\kappa =0, \nu_\kappa =0$ and $g_\kappa=0$, and from the assumptions, $0\preceq P_k^{(\kappa)}\uparrow L, k=0,\ldots,\kappa$ as $\kappa \to \infty$, for $L$ that solves $(I-\Ric^\alpha)(U)=0$. Therefore, 
\begin{equation}\label{eq:limit.norm.alpha}
  \lim_{\kappa\to\infty} E\bigl[\sum_{k=0}^{\kappa-1}\alpha^k\|y_k\|^2|x_0\bigr]=
 \\V^{u^*}_0(x_0) =\|x_0\|_L+\langle \mu_0^{(\infty)}, x_0\rangle + g_0^{(\infty)}   
\end{equation}
 and thus, 
\[
    \mathcal{E}_2^{\alpha*} = g_0^{(\infty)} = E\Bigl[\sum_{k=0}^\infty \alpha^k( \alpha\varpi_1(L)  +  \rho_k) \big|x_0\Bigr] 
\]    
 
 For the  $H_2$-nom in C\`esaro mean form, the assumptions in the theorem hold with $\alpha=1$, which implies that they also hold to each $\alpha<\bar{\alpha}$ for some $\bar{\alpha}>1$.  In addition, since  $k\to \rho_k$ (or, explicitly, $k\to\rho^\alpha_k$) are bounded processes for each $\alpha< \bar{\alpha}$,
 one can apply Theorem 3.6 of \cite{CSVIU:norm} with $Q=C^\trp C$ to provide the vanishing discount formula,
 \begin{multline*}\allowdisplaybreaks
\mathcal{E}_2^{*} =\lim_{\kappa\to \infty} \frac{1}{\kappa} \sum_{k=0}^{\kappa-1} E[\|y^*_k\|^2]=
\lim_{\alpha \to 1}(1-\alpha)\mathcal{E}_2^{\alpha*}= 
\\
\lim_{\alpha \to 1}(1-\alpha) E\Bigl[\sum_{k=0}^\infty \alpha^k\left( \alpha\varpi_1(L)  +\rho^\alpha_k\right)|x_0=0\Bigr]
=\varpi_1(L) + \lim_{k \to \infty}E[\rho^1_k] 
\end{multline*}
 which shows the expressions of the optimal norm in \eqref{eq:optomal.norms}.

 For the overtaking criteria, note first that if one sets $P_\kappa = L$, the solution of $(I-\Ric^\alpha)(U)=0$, and $\mu_\kappa = \mu_\kappa^{(\infty)}$ as in \eqref{eq:vk.infty}, it follows that $P_k=L$, $\mu_k = \mu_\kappa^{(\infty)}$ hold for each $k\le\kappa$. In this situation, we obtain the expressions of  $k\to \rho_k, \bar{u}_k$ exactly as in \eqref{eq:rho.L} and \eqref{eq:argument.optimization}, respectively. 
 
 On the other hand, $\{u^*_k\}_{k=0,\ldots,\kappa}$ satisfying \eqref{eq:optimal.law}, with $\varrho_k(u)$ defined by \eqref{eq:rho.u}, attains the infimum of the finite horizon criteria $\mathcal{E}_2^{\alpha,\kappa}(y(\cdot))$ for each $\kappa$,  as  \eqref{eq:optimal.finite.cost.comp} and \eqref{eq:optimalcost.kappa} indicate, and we can set the comparision $V_0^{u^*} \le \mathcal{E}^{\alpha,\kappa}_2(\bar{y}(\cdot))$.
 In this connection, we have set  $P_\kappa = 0$ and $\mu_\kappa =0$, and we know from part (I) that 
 $P^{(\kappa)}_k\uparrow L$ and $\mu^{(\kappa)}_k \to \mu^{(\infty)}_k$ as $\kappa\to\infty$, with $L$ the solution of $(I-\Ric^\alpha)(U)=0$.
Convergence leads to  convergence of $\varrho_{k}(u)$ in \eqref{eq:rho.u} to $\rho_{k}(u)$ in \eqref{eq:rho.L} and hence, $\varrho^*_k \to\rho_k$ in \eqref{eq:stage.optimization}.  One can conclude that for any $\epsilon>0$, there is a sufficiently large number $\kappa_0$ such that we can compare $V_0^{u^*} > \mathcal{E}^{\alpha,\kappa}_2(\bar{y}(\cdot))-\epsilon$ for any $\kappa\ge \kappa_0$.
\end{proof}

 \begin{proof}[Proof of part (III)]
 For this proof we rely on the fact that $k\to\rho_k$ is a bounded process and we can read the results of the norm of system $\Theta$ in \cite{CSVIU:norm} by simply replacing $ \varpi(L)$ therein by $ \varpi_1(L) + \varrho_k$ for the accumulated value $k\to g_k$ in \eqref{eq:stage.residual} at each stage.  Moreover, note that  $L\succeq0$ the solution of $(I-\Ric^\alpha)(U)=0$ is also the solution of $(I-\mathcal{H}^\alpha_G)(U)$ for $G=- \Lambda(L)^{-1} \Sigma(L)$, which implies that $\alpha^{1/2}(A+BG)$ is d-stable relative to $\alpha(\mathcal{Z}_x+G^\trp\mathcal{Z}_uG)$, see Appendix \ref{app:stabilization.G}. 
 
 Let us consider first $\alpha<1$ and \eqref{eq:limit.norm.alpha} implies, using the notation of $Q$-mean energy in \cite{CSVIU:norm} that $\lim _{k \rightarrow \infty} \mathcal{E}_{2, Q}^{\alpha, k}(x(\cdot))<\infty$  for $Q={C}^{\trp} {C}$ and any $x(0)$. Then \cite[cor 2.6 ii]{CSVIU:norm} applies to guarantee $\alpha$-stochastic stability. When $\alpha\ge1$ one can set $P_\kappa = L$ and $\mu_\kappa = \mu^{(\infty)}_\kappa$ and write that 
 \begin{multline}
   E\bigl[\sum_{k=0}^{\kappa-1}\alpha^k\|y_k\|^2|x_0\bigr]\le
 \|x_0\|_L+\langle \bar{\mu}, |x_0|\rangle + \sum_ {k=0}^{\kappa-1} \alpha^{k+1} ( \varpi_1(L) + \varrho_k)
 \\
 \le c_0\|x_0-\zeta\|^2 + \kappa \alpha^{\kappa} c_1
\end{multline}
 with $\zeta=-\frac{1}{2}(\tilde{L})^{-1}\bar{\mu}\sqcdot\mathcal{S}(x_0)$, $\tilde{L}=L+\epsilon I$, $c_0 =\lambda^+(L)$ and $c_1 = \alpha(\varpi_1(L) + \max_k|\varrho_k|)$. Hence, $\alpha$-stability follows when $\alpha\ge1$. Stochastic stability for the C\`{e}saro mean follows simply from the above with $\alpha=1$.
  \end{proof}

\begin{remark}
The positive semidefinite solution of the modified Riccati equation together with the assumption on the spectral radius when $\alpha>1$ in Theorem \ref{theo:controlled.solution} unveil a solution that relies on a form that resembles the deterministic Riccati solution, cf.  \eqref{eq:rho.L}.  However, it is stochastic in nature due to the dependence on the signals $k \to \mathcal{S}(x_k), \mathcal{S}(u_k)$ along the path. Note that if $\bar{\sigma}_x$ and $\bar{\sigma}_u$ are null, one retrieves the discounted and the long-run $H_2$-norm problems in one hand, or the optimal overtaking problem of a linear system driven by white noise on the other.

We can also retrieve from the general scenario drawn by Theorem \ref{theo:controlled.solution} the idea of asymptotic solutions developed in \cite{bib:3}. The first aspect is that for each control action $u_i$ in $u=[\begin{smallmatrix}u_1& \cdots& u_m \end{smallmatrix}]^\trp$ the regions $\mathcal{R}_i^+$ and $\mathcal{R}_i^-$ are separated by an inaction region $\mathcal{R}_i^0$. It will become clear in the next section that $\mathcal{R}_i^0$ may form a $n$-dimensional region in the state space, $\mathds{R}^n$. 

Whenever the state evolves in some regions for which all the signs are invariant during a sufficiently large number of time stages, say, $k \to \mathcal{S}(x_k)\simeq \bar{s}_x$ and $k \to \mathcal{S}(u_k)\simeq \bar{s}_u$, one can try to approximate, 
\begin{multline*}
 \mu_0 =  E\Bigl[\sum_{k=0}^\infty \alpha^{k+1} \mathcal{A}^k\bigl(\mathcal{W}_{x}(L) \mathcal{S}(x_k) 
+G^\trp\mathcal{W}_{u}(L)\mathcal{S}(u^*_k)\bigr) |x_0\Bigr]
    \\
    \simeq  {\mu}(\bar{s}_x,\bar{s}_u):=\alpha(I-\mathcal{A})^{-1}\bigl(\mathcal{W}_{x}(L) \bar{s}_x +G^\trp\mathcal{W}_{u}(L)\bar{s}_u\bigr)
\end{multline*}
 with $G=-\Lambda(L)^{-1}\Sigma(L)$ and $\mathcal{A}= A+BG$. In this situation, the optimal feedback solution can simply be approximate by
 \begin{equation}\label{eq:assymp.feedback}
  u^0(\bar{s}_x,\bar{s}_u)(x):= 
  \\
  - \Lambda(L)^{-1}
\left( \Sigma(L)x+ \frac{1}{2}(B^\trp {\mu}(\bar{s}_x,\bar{s}_u) + \mathcal{W}_{u}(L) \mathcal{S}(u^0))\right)
 \end{equation}

From that comes the idea of asymptotic solutions of \cite{bib:3} and \cite{bib:4}. First, attach appropriate signal vectors $\bar{s}_u=[\begin{smallmatrix}
\bar{s}_{u,1}& \cdots & \bar{s}_{u,m}
\end{smallmatrix}]^\trp$ with each $\bar{s}_{u,i}\not=0$ and consider the regions  $\mathcal{R}^{\bar{s}_{u}} =\cap_i\mathcal{R}_i^{\bar{s}_{u,i}}$ of $\mathds{R}^n$. Now, inside the intersection of such regions $\mathcal{R}^{\bar{s}_{u}}$  with each octant of $\mathds{R}^n$, the feedback solution in \eqref{eq:assymp.feedback} is \emph{asymptotically optimal} as $\|x\|\to \infty$ and far away from the boundaries of those regions. 

Finally, note that when each signal in $\bar{s}_x$ and $\bar{s}_u$ is nonzero, the true value $\mu_0$ belongs to the intervals created by  $\mu(\bar{s}_x,\bar{s}_u)$ in one extreme, and the zero vector on the other. Thus, the vectors ${\mu}(\bar{s}_x,\bar{s}_u)$ form upper values (in magnitude terms) for the true value of $\mu_0$. 
\end{remark}

\section{Nondifferentiable optimization and the inaction region}\label{se:inaction.representation}

The inaction regions $\mathcal{R}^0 = \cap_{i=1}^m \mathcal{R}_i^0$ can be characterized with the elements of Theorem \ref{theo:controlled.solution}. Suppose that $\mathcal{S}(u)\equiv 0$ is optimal. Then, from \eqref{eq:rho.min},
\begin{equation*}
     \rho_k = \|u^0_k\|^2_{\Lambda(L)} 
     -\frac{\alpha}{4}\left\|B^\trp \mu_{k+1}  \right\|^2_{\Lambda(L)^{-1}} 
\end{equation*}
is readily satisfied. Clearly, each inaction region $\mathcal{R}_i^0$ can be characterized by setting  $u_i=0$ and assume that $\rho_k$ is attained in such a way. However, that allows us a glimpse on the inaction region, but not a way of determining it. Note also that it is hard to tell if the problem of minimizing  \eqref{eq:stage.optimization} is amenable or not.
 
For the purpose of having a better view on the optimal solution we first quote \eqref{eq:dif.one.step.expected} in Lemma \ref{lemm:variation.Vu}, which is indeed the base for the proof of Theorem \ref{theo:controlled.solution}. 
To seek the optimal control solution at the $k$-th stage is precisely the problem of minimizing the difference in \eqref{eq:dif.one.step.expected} wrt the choice of $u_k=u$.   We  retrieve the absolute values involved in the original form, to write for a $u_k=u$, the difference $J_u:\mathds{R}^m\to \mathds{R}$, as
\begin{multline}\label{eq:optimal.l1l2}\allowdisplaybreaks
     J_u := \alpha^{-k} E\left( V^u_{k+1}( x_{k+1}) - V^u_k (x_k) +\|y_k\|^2 |x_k\right)=
\\
      \|x_k\|^2_{M} + \langle  m^1_k, x_k \rangle   + \langle  m^2_k, |x_k| \rangle 
        + \|u\|^2_{\Lambda(L)} +\alpha \langle B^\trp \mu_{k+1}  +2\Sigma(L) x_k,u \rangle 
    \\
    +\langle \alpha \mathcal{W}_{u_d}(L)  ,|u| \rangle 
         +\alpha g_{k+1} + \alpha \varpi_{\textcolor{blue}{1}}(L) - g_{k}
         \\
         =  \|u\|^2_{\Lambda(L)} +\alpha\langle B^\trp \mu_{k+1}  +2\Sigma(L) x_k,u \rangle +
         \langle \alpha \mathcal{W}_{u_d}(L)  ,|u| \rangle + f_k
\end{multline}
where $M, m^1_k,m^2_k,f_k$ denote other terms in \eqref{eq:dif.one.step.expected} that does not depend on $u$.
For the sake of $u\to J_u$ being a strictly convex function, we pinpoint the relevant assumptions. The first is a standing assumption in Section \ref{sec:optimal.solution}, recall that $D^\trp D\succ0$ is required in Theorem \ref{theo:controlled.solution}.
\begin{assumptions}\label{assump:convex.Ju} 
     $\Lambda$ is invertible and   $\mathcal{W}_{u}$ is a linear-positive operator.
\end{assumptions}

 Note that the assumption on  $\mathcal{W}_{u}$ implies that the $\mathcal{W}_{u_d}(U)\in\mathds{R}^m$ has nonnegative entries for each $U\succeq0$.
With these requirements, it is easy to check that  $u\to J_u$ is convex but nondiferentiable and the optimality condition can be expressed as,
\begin{equation}\label{eq:sing}
0\in\partial J(u)
\end{equation}
where $\partial J(u)$ stands for the Clarke's generalized gradient set, defined for a locally Lipschitz function or equivalently, the Rockafellar's subgradient set, when $J$ is convex function, cf. \cite{Clarke,Rockafellar}. 
From Lemma 3.2 of \cite{bib:3}, 
\begin{equation*}
\mathcal{R}_i^0=\bigl\{x\in\mathds{R}^n:\lim\limits_{u_i\uparrow0}\partial_{u_i}J(u)\leq 0\leq\lim\limits_{u_i\downarrow0}\partial_{u_i}J(u)\bigr\}
\end{equation*}
where $\partial_{u_i} J(u)$ indicates Clarke's gradient of the $i$-th component of the multidimensional function $J(u)$. 
To get the minimum of $J_u$, we recast \eqref{eq:sing} as a special type of normal equation. Denote in \eqref{eq:optimal.l1l2},
\begin{equation}\label{eq:simpler.elements}
W:=\frac{1}{2}\Lambda(L)^{-1},\quad b:= \alpha( B^\trp \mu_{k+1}  +2\Sigma(L) x_k)^\trp, \quad c:= \alpha\mathcal{W}_{u_d}(L)
\end{equation}
and we can express,
$
\partial J(u) = W^{-1} u + b + c\sqcdot \partial(|u|)
$, 
where  $\partial(|u|)$ represents the set of the generalized gradient of the absolute values involved in \eqref{eq:optimal.l1l2}. Namely, for $u\in \mathds{R}^m$, $\partial(|u|)= [\begin{smallmatrix}\partial(|u_1|) & \cdots & \partial(|u_m|)\end{smallmatrix}]^\trp$ is such that 
\[
\begin{cases}
\partial(|u_i|) =\sign(u_i),&\text{ if $u_i\not= 0$},
 \\
\partial(|u_i|)  \in [-1,+1],& \text{ if $u_i= 0$}.
\end{cases}
 \qquad i=1,\ldots,m.
 \]

To express the optimality condition in \eqref{eq:sing} we consider the \emph{generalized normal equation}
\begin{equation}\label{eq:gen.normal}
 W^{-1} u + b + \gamma=0, 
\end{equation}
 for some $\gamma\in\{c\sqcdot \xi: \xi\in \partial(|u|) \}$, with the following understanding, 
 \[
 \begin{cases}
\gamma_i = c_i \sign(u_i),&\text{ if $u_i\not= 0$},
 \\
 \gamma_i  \in (-c_i,+c_i),& \text{ if $u_i= 0$}.
\end{cases}
\quad i=1,\ldots,m,
\]
in which we set a one-to-one correspondence at the extremes $\pm c_i$ for each $i$. One can write \eqref{eq:gen.normal} equivalently as
$ (W^{-1} \diag(\theta)  +I)\gamma+ b=0$, for $\theta=[\begin{smallmatrix} \theta_1 & \cdots& \theta_m\end{smallmatrix}]^\trp$, where
 \[
 \theta_i = \begin{cases}
 {|u_i|}/{c_i}\ge 0,&\text{if $c_i\not=0$,}
 \\
 0,&\text{otherwise.}
 \end{cases} \qquad i=1,\ldots,m.
 \]
Hence, \eqref{eq:gen.normal} is equivalent to 
 \begin{equation}   \label{eq:normal2}
    ( \diag(\theta)  +W)\gamma+ Wb=0, 
 \end{equation}

Eq. \eqref{eq:normal2} is an equivalent form of the generalized normal equation for the mixed $\ell_2$/$\ell_1$ problem that arises in connection with the cost $J_u$. Note however, that contrary to  the case of standard normal equation it does involve more than a simple matrix inversion.  
A special procedure to attain the minimum of $J_u$ is devised here, which resembles SOR iterations \cite{Hackbusch}. In what follows, we explicitly use the fact that $\theta$ is a vector with nonnegative entries.

Let us refer to the $(i,j)$-th element of matrix $U$ as $U_{i,j}$, and recall the notation $U_d=\Diag(U)\in\mathds{R}^m$.  We set $\min\{\cdot,\cdot\}$ and $\max\{\cdot,\cdot\}$ operations of two vectors of same dimension as the componentwise corresponding operations.
\begin{theorem}\label{th:matrix.inverted}
For any $z^0\in \mathds{R}^m$ and some $0<\omega<2$, consider the recursion
\begin{equation}\label{eq:main.recursion.1}
z^k_i=  \left(1-\omega\right)z^{k-1}_i
\\
-\frac{\omega}{W_{i,i}} \Bigl(\sum_{j=1}^{i-1} W_{i,j}\left(\gamma_j^{k}+b_j\right)  + \sum_{j=i+1}^{n} W_{i,j}\left(\gamma_j^{k-1} + b_j\right)\Bigr) -\omega b_i
\end{equation}
with
\begin{align}\label{eq:gamma.x}
\gamma^{k} &= \min\left\{ c, |z^k|\right\}\sqcdot \sign(z^k),
\\
\label{eq:v.x}
\nu^k & = \max\left\{0, W_d\sqcdot(|z^k| -c)\right\}\sqcdot \sign(z^k) 
\end{align}

Then $\gamma^k ,\nu^k\to \gamma^*, v^*$ that satisfy the generalized normal equation \eqref{eq:gen.normal}.

\end{theorem}

\begin{proof}
Note that \eqref{eq:normal2} can be written equivalently as,
\begin{multline}\label{eq:consistency.recursive2}
0 = (\diag(\theta)  +W)\gamma+ Wb = 
\\
\frac{1}{\omega}(\diag(\theta)  +{W}_d)\gamma +
\Bigl(1- \frac{1}{\omega}\Bigr)(\diag(\theta)  +{W}_d)\gamma 
+ ({W}_U+{W}_L)\gamma  + W b
\end{multline}
for $\omega \not=0$, where  ${W}_U$ and ${W}_L$ are the strictly upper and lower  triangular matrices respectively originated from $W$. 
Now, one employs \eqref{eq:consistency.recursive2} to write component-wisely the relaxation for increasing $i=1,\ldots,m$, 
\begin{equation}\label{eqs:alg.no.inv}
\Bigl(\frac{1}{W_{i,i}} \theta^{k}_i +1\Bigr)\gamma^{k}_i =
\\
\left(1-\omega\right)\Bigl(\frac{1}{W _{i,i}} \theta^{k-1}_i +1\Bigr)\gamma^{k-1}_i + \omega\ell^k_i
\end{equation}
for $\omega>0$, where,
\begin{equation}\label{eq:rec.gamma}
\ell^k_i= -\frac{1}{W_{i,i}} \Bigl(\sum_{j=1}^{i-1} {W}_L(i,j)\left(\gamma_j^{k}+b_j\right) 
+ \sum_{j=i+1}^{n} {W}_U(i,j)\bigl(\gamma_j^{k-1} + b_j\bigr)\Bigr) -b_i
\end{equation}

For simplicity let us represent the rhs of \eqref{eqs:alg.no.inv} as 
\[
\ell^k_i(\omega)= \left(1-\omega\right)\Bigl(\frac{1}{W _{i,i}}\theta^{k-1}_i +1\Bigr)\gamma^{k-1}_i + \omega\ell^k_i
\]
Suppose that $c_i>0$, and note that if $|\ell^k_i(\omega)|<c_i$, then $\gamma_i^k = \ell_i^k(\omega), \theta_i^{k} = 0$ and therefore, $\nu_i^{k}: =  \theta_i^{k}\gamma_i^k=0 $. Otherwise, 
\begin{equation}\label{eq:no-gap}
\left(\frac{1}{W _{i,i}} \theta^{k}_i +1\right)\gamma^{k}_i =\ell^k_i(\omega)
\end{equation}
with $\gamma_i^{k} = c_i \sign(\ell^{k}_i(\omega))$ and $\theta^{k}_i  = \frac{1}{\gamma^{k}_i} W_{i,i}(\ell^k_i(\omega)- \gamma^{k}_i)$.
When $c_i=0$, \eqref{eq:no-gap} imposes that $v_i^k := W_{i,i}(\ell^k_i(\omega)- \gamma^{k}_i)$ with $\gamma^{k}_i=0$. 

For the moment, assume that $c_i>0, \forall i$. The above relations read as
\begin{align*}
\gamma^{k} &= \min\bigl\{ c, |\ell^k(\omega)|\bigr\}\sqcdot \sign(\ell^k(\omega))
\\
\theta^{k} &= \max\bigl\{ 0, {W}_d\sqcdot\bigl(\diag(c)^{-1}|\ell^k(\omega)| - \mathds{1}\bigr) \bigr\}\sqcdot \sign(\ell^k(\omega))
\end{align*}
 Here, $\mathds{1}$ represents the $m$-dimensional vector made up by one at each component. Denote $z^k:=\ell^k(\omega)$
and the recursion above is determined precisely as in \eqref{eq:main.recursion.1} with \eqref{eq:gamma.x}.
Set $\nu^k := \theta^k \sqcdot \gamma^k$ and note  that  $\nu^k$ can be expressed as in \eqref{eq:v.x}.

To complete this part of the proof, it remains to verify the case when $c_i=0$. But notice that \eqref{eq:gamma.x} and \eqref{eq:v.x} yield $\gamma_i^k=0$ and $\nu^k_i  = W_{i,i}|z_i^k|\sign(z_i^k)$ as required.

\smallskip

For the proof of convergence,  we mention that the devised algorithm resembles SOR iterations, and the standard SOR method is brought here to this end.
 It can be written component-wisely  similarly to \eqref{eq:main.recursion.1}, as
\begin{equation*}
   x_i^k = (1-\omega) x^{k-1}_i -\frac{\omega}{W_{i,i}}\Bigl( \sum_{j=1}^{i-1} W_{i,j}\bigl(x_j^{k}+b_j\bigr) + \sum_{j=i+1}^{n} W_{i,j}\bigl(x_j^{k-1} +b_j\bigr)    \Bigr)-\omega b_i,
\end{equation*}
for $i=1,\ldots,m$ and $W\succ0$. The recursion converges monotonically in the energy (or spectral) norm for $0<\omega<2$ for any $x^0$, cf. \cite{Hackbusch}. The particular case when $\omega=1$ corresponds to the  Gauss-Seidel method\footnote{These methods are normally not employed in the present form, since their main  motivation is  to avoid the inversion of matrix $2\Lambda(L) = W^{-1}$ in \eqref{eq:simpler.elements} for large problems, hence beyond the scope here.}.  

Let us rephrase the SOR iterations method by defining the linear operator $\mathcal{T}^\omega:\mathds{R}^m\to\mathds{R}^m$, such that 
\[
\mathcal{T}^\omega x = (I+ \omega \diag({W}_d)^{-1} {W}_L)^{-1}(I-\omega(I + \diag({W}_d)^{-1}{W}_U))x,
 \]
 for $x\in\mathds{R}^m$,
 in such a way that we can represent the SOR iterations as
\[
x^k = \mathcal{T}^\omega(x^{k-1}) -  \omega(\diag({W}_d)+ \omega  {W}_L)^{-1}Wb
\]
Convergence for the SOR comes from the fact that the spectral matrix norm $ \|\mathcal{T}^\omega\|<1$ for $0<\omega<2$, which implies that $\mathcal{T}^\omega$ is a strict contraction mapping operator in the norm above. The sequence $\{x^k\}$ tends  to the unique fixed point of the recursion 
$
\bar x = - {\omega}( I -\mathcal{T}^\omega)^{-1} (\diag({W}_d) + \omega  {W}_L)^{-1}  b.
$

To express the actual method in operator form, consider the nonlinear operator $\vartheta:\mathds{R}^m\to\mathds{R}^m$, such that for $x\in\mathds{R}^m$,
\[
\vartheta(x) = \min\{c,|x|\}\sqcdot\sign(x)
\]

\begin{figure}[!h]
\begin{center}

\tikzset{every picture/.style={line width=0.75pt}} 

\begin{tikzpicture}[x=0.75pt,y=0.75pt,yscale=-0.8,xscale=0.8]

\draw    (124.25,159.75) -- (454.25,159.25) ;
\draw [shift={(456.25,159.25)}, rotate = 539.9100000000001] [color={rgb, 255:red, 0; green, 0; blue, 0 }  ][line width=0.75]    (10.93,-3.29) .. controls (6.95,-1.4) and (3.31,-0.3) .. (0,0) .. controls (3.31,0.3) and (6.95,1.4) .. (10.93,3.29)   ;
\draw    (285.25,230.25) -- (285.25,78.25) ;
\draw [shift={(285.25,76.25)}, rotate = 450] [color={rgb, 255:red, 0; green, 0; blue, 0 }  ][line width=0.75]    (10.93,-3.29) .. controls (6.95,-1.4) and (3.31,-0.3) .. (0,0) .. controls (3.31,0.3) and (6.95,1.4) .. (10.93,3.29)   ;
\draw [color={rgb, 255:red, 208; green, 2; blue, 27 }  ,draw opacity=1 ][line width=1.5]    (360.25,100.25) -- (208.75,220.25) ;
\draw    (283,220.5) -- (287.25,220.75) ;
\draw    (283,100.5) -- (287.25,100.75) ;
\draw [color={rgb, 255:red, 208; green, 2; blue, 27 }  ,draw opacity=1 ][line width=1.5]    (126.25,220.75) -- (208.75,220.25) ;
\draw [color={rgb, 255:red, 208; green, 2; blue, 27 }  ,draw opacity=1 ][line width=1.5]    (360.25,100.25) -- (442.75,99.75) ;
\draw  [dash pattern={on 0.84pt off 2.51pt}]  (391.75,76.25) -- (191.75,233.25) ;

\draw (253.5,94.5) node [anchor=north west][inner sep=0.75pt]  [font=\footnotesize]  {$+c_{i}$};
\draw (253,214) node [anchor=north west][inner sep=0.75pt]  [font=\footnotesize]  {$-c_{i}$};
\draw (439,165.5) node [anchor=north west][inner sep=0.75pt]  [font=\footnotesize]  {$x$};
\draw (290.5,75) node [anchor=north west][inner sep=0.75pt]  [font=\footnotesize]  {$\vartheta_i ( x)$};

\end{tikzpicture}

\caption{$\vartheta_i(x) = \min\{c_i, |x_i| \}\sign(x_i),  i=1,\ldots,m$. }
\label{fig:phi.estim2}
\end{center}
\end{figure}
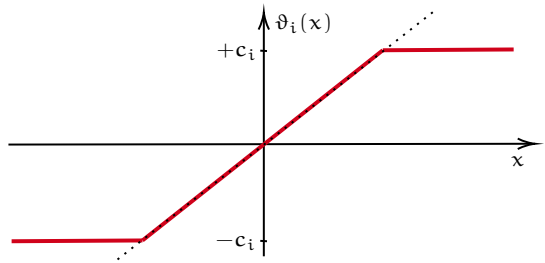

From the fact that $|\vartheta(x)| =  \min\{c,|x|\} \le |x|$ and $\|\vartheta(x) - \vartheta(y)\|\le \|x-y\|$ one has that $\vartheta$ is a contractive mapping operator, see Fig. \ref{fig:phi.estim2}. Note that the recurrence proposed in theorem can be expressed as
\begin{equation}\label{eq:recurr.2}
\gamma^k = \vartheta\circ \Bigl(\mathcal{T}^\omega(\gamma^{k-1}) - \omega(\diag({W}_d)+ \omega  {W}_L)^{-1}Wb\Bigr)
\end{equation}
and since $\mathcal{T}$ is  a strictly contractive operator, it holds for the compound operator $\vartheta\circ\mathcal{T}^\omega$ that
$\|\vartheta\circ\mathcal{T}^\omega\| \le \|\vartheta\|\|\mathcal{T}^\omega\|$ with  $\|\vartheta\|\le1$ and $\|\mathcal{T}^\omega\|<\beta^\omega$ for some $0<\beta^\omega<1$. 

Hence $\vartheta\circ\mathcal{T}^\omega$ is  strictly contractive mapping in the energy norm. As a result there is an unique fixed point for the recurrence in \eqref{eq:recurr.2} that satisfies,
$$
\gamma^* = \vartheta\circ\Bigl(\mathcal{T}^\omega(\gamma^{*})-\omega(\diag({W}_d)+ \omega  {W}_L)^{-1}Wb\Bigr)
$$ 
and for any $\gamma^0\in \mathds{R}^m$, $\gamma^k\to \gamma^*$ in the euclidean norm. Since $\{\gamma^k\}$ converges, 
\[
z^k = \mathcal{T}^\omega(\gamma^k)-\omega(\diag({W}_d)+ \omega  {W}_L)^{-1}Wb
\]
also necessarily converges. 
  Now, from \eqref{eq:v.x} it comes that the sequence $\nu^k\to v^*$, the unique solution of the generalized normal equation 
 $$
 ( \diag(\theta)  +W)\gamma+ Wb=  \nu  +W(\gamma+ b)=0.
 $$
 with $-c_i\le \gamma_i\le +c_i$.
 \end{proof}

We are now in position to refer directly to the original optimal control problem. Let $b_i$ indicates the  $i$-th element of vector $b$ and $A_i$ the $i$-th row of matrix $A$.

\begin{theorem} \label{th:solution}
Under the Assumptions \ref{assump:convex.Ju} and in Theorem \ref{theo:controlled.solution} 
the accumulation point $v^*$ of the recursion in Theorem \ref{th:matrix.inverted} is such that $J_{v^*}\le J_u, \forall u\in\mathds{R}^m$. Moreover,  
\begin{equation}\label{eq:generalized.normal}
 v^*=-\Lambda(L)^{-1}\Bigl( \Sigma(L) x_k +\frac{1}{2} \left(B^\trp \mu_{k+1} +  \gamma^*\right)\Bigr)
\end{equation}
for  $\gamma^*=[\begin{smallmatrix}
\gamma_1^*  & \cdots & \gamma_m^*
\end{smallmatrix}]^\trp\in\mathds{R}^m$ such that $|\gamma^*|\le \mathcal{W}_{u_d}(L)$ understood componentwisely. 
In addition,  $\bar{u}_k=v^*$ is the optimal control in the sense of Theorem \ref{theo:controlled.solution}. 

Provided that $|\gamma_i^*|< \left(\mathcal{W}_{u_d}(L) \right)_i$ then $\bar{u}_i =0$ and each $x\in \mathcal{R}_i^0$ satisfies
\[
 \bigl|2\,\langle \Sigma(L)_i ,x\rangle+ \langle B_i, \mu_{k+1}\rangle\bigr| <   \left(\mathcal{W}_{u_d}(L) \right)_i
\]

\end{theorem}

\begin{proof}
		The result comes from the correspondence made in \eqref{eq:simpler.elements}, the generalized normal equation and the recursion developed in Theorem \ref{th:matrix.inverted}. The expression for a point belonging to region $\mathcal{R}_i^0$ is obtained by setting $v^*_i=\bar{u}_i=0$ in \eqref{eq:generalized.normal} and  noting that $\gamma_i^*\in(-c_i,+c_i)$ necessarily, with $c_i=\left(\mathcal{W}_{u_d}(L) \right)_i$.
\end{proof}

\begin{remark} 
The boundaries of the region $\mathcal{R}_i^0$ indicated in Theorem \ref{th:solution} seems to be confined by parallel hyperplanes. Note however, that $\mu_{k+1}=E[v_{k+1}|x]$ is determined  by \eqref{eq:reccurrence.mu}, which depends upon the very state point $x_k=x$. 

We should emphasize that the generalized normal equation in \eqref{eq:generalized.normal} is more adequate to evaluate the optimal control than the direct approach to the problem in \eqref{eq:rho.min}. The type of nondifferentiable coordinate descent method presented in Theorem \ref{th:matrix.inverted} can efficiently tackle the solution.

\end{remark}

\section{Conclusion}
\label{sec:conclusion}

The paper brings the fundamentals of CSVIU modeling as a contribution to the uncertain systems literature and robust theory. The focus is the optimal control problem, studied in terms of energy and power $H_2$-norms and the optimal overtaking criterion. It characterizes the optimal stochastic stabilizing solution in the proposed senses. A detectability notion and a perturbed Riccati equation solution play a fundamental role here. However, for the power norm and the counter-discounted problems, a condition on the spectral radius of a closed-loop matrix is further required. 

The optimal solution is global and not seen in previous results. Theorem \ref{theo:controlled.solution} reveals the implication of a perturbed Riccati-like feedback form in the local optimization solution. In a second standing, Theorem \ref{th:solution} resorts to nondifferentiable tools in the corresponding problem, taking the Riccati-like feedback form into a different understanding. The underlying nondifferentiable normal equation exposes the rise of the inaction region of control, and it can be solved, leading to a procedure to attain the optimal solution. 

These findings deliver a complete picture of the $H_2$ energy and power norms, overtaking control problems, and the relation between counter-discounted, discounted, and long-run formulations. They also equip us to solve the optimal control problem for the CSVIU model $\Theta_{\textrm{ctr}}$.

\bibliography{Daniel}

\newpage
\appendix

\section{Proof of Lemma \ref{lemm:variation.Vu} }\label{app:lemm.variation.Vu}

We retrieve similar evaluations of those in the proof of   \cite[Lem\;2.1]{CSVIU:norm}, driven by an exogenous process. 
    Having in mind that  $\langle r, |x| \rangle = \langle \mathcal{S}(x), r \sqcdot x \rangle= \langle\mathcal{S}(x)\sqcdot r,  x \rangle$ and  $\tr \{ U  \diag(|x|)\}= \langle \mathcal{S}(x), \Diag(U) x \rangle$ hold true for any $U \in \mathds{S}^{n+}$ and $r,x \in \mathds{R}^n$, and accounting for  the dynamics of system $\Theta_{\textrm{ctr}}$, one can evaluate the variation of $V^u$ in \eqref{eq:Vu} along a path $k\to x_k$ (with the compact notation in \eqref{eqs:compact.theta}),
\begin{multline*} \allowdisplaybreaks 
          \alpha^{-k} \big(  V^u_{k+1}( x_{k+1})- V^u_k (x_k) \big)=
          \\ 
          \alpha \| x_{k+1}\|^2_{P_{k+1}} +  \alpha \langle s_{k+1}, r_{k+1} \sqcdot x_{k+1} \rangle  + \alpha g_{k+1} -
          (\| x_{k}\|^2_{P_{k}} + \langle s_k, r_{k} \sqcdot x_{k} \rangle + g_k )=
\\
\alpha \| A x_k + B u_k \|^2_{P_{k+1}}   +  2 \alpha (A x_k +B u_k)^\trp P_{k+1} \sigma(x_k) \omega_{0_k}
           + \alpha \| \sigma(x_k)\omega_{0_{k}} \|^2_{ P_{k+1} }  - \| x_k \|^2_{P_k } +
          \\
      \alpha \langle s_{k+1}, r_{k+1} \sqcdot \big( A x_k + B u_k +  \sigma(x_k)\omega_{0_k} \big) \rangle - \langle s_k, r_{k} \sqcdot x_{k}  \rangle +  
          \alpha g_{k+1} - g_{k} 
\end{multline*}
where we set $s_{k} = \mathcal{S}(x_{k})$ and $s_{k+1} = \mathcal{S}(x_{k+1})$. Note that
\begin{multline*}
E  \bigl[\| \sigma_x(x_k)\omega_{0_k} + \sigma_u(u_k)\epsilon^u_k \|^2_U|x_k =x, u_k=u\bigr]  =
     \| x \|^2_{\mathcal{Z}_x(U)} + \langle \mathcal{S}(x), \mathcal{W}_x(U)x \rangle  
\\
+\| u \|^2_{\mathcal{Z}_u(U)} + \langle \mathcal{S}(u), \mathcal{W}_u(U)u \rangle + \varpi_1(U)
\end{multline*}
holds. Set also $s_u:=\mathcal{S}(u_k)$ and evaluate,
\begin{multline*}\allowdisplaybreaks
      \alpha^{-k} \big( V^u_{k+1}( x_{k+1})- V^u_k (x_k)\big) +\|y_k\|^2 =
      \\
     \alpha \big( \| A x_k \|^2_{P_{k+1}} + \| x_k \|^2_{\mathcal{Z}_x(P_{k+1})}\big) + \| C x_k \|^2 +
      \alpha \big( \| B u_k \|^2_{P_{k+1}} + \| u_k \|^2_{\mathcal{Z}_u(P_{k+1})}\big) +  \| D u_k \|^2 +
     \\
        \alpha \langle s_{k+1}, r_{k+1} \sqcdot (A x_k+ Bu_k) \rangle + \alpha\langle   \mathcal{W}_{xd}(P_{k+1}) \sqcdot s_k,  x_{k} \rangle +
    \\
        \alpha\langle   \mathcal{W}_{ud}(P_{k+1}) \sqcdot s_u,  u_{k} \rangle + 2\langle D^\trp C x_k,u_k\rangle
    \\
    \alpha g_{k+1} + \alpha \varpi_1(P_{k+1})       - \| x_k \|^2_{P_{k} } -  \langle s_{k}\sqcdot r_{k} , x_{k} \rangle - g_{k} + m_k =
 \end{multline*}
 \vspace{-5.3ex}
\begin{multline*}
      x_k^\trp \big(\mathcal{L}^\alpha(P_{k+1})+ C^\trp C -  P_{k} \big) x_k + \alpha u_k^\trp \Lambda(P_{k+1}) u_k+
     \\
        \langle \alpha A^\trp r_{k+1} \sqcdot s_{k+1} + \alpha \mathcal{W}_{x}(P_{k+1})s_k - s_k\sqcdot r_{k}   ,x_{k} \rangle + 
    \\
       \alpha \langle  B^\trp r_{k+1} \sqcdot s_{k+1} + \mathcal{W}_{u}(P_{k+1})s_u +2\Sigma(P_{k+1})x_k ,u_{k} \rangle + 
    \\
         \alpha g_{k+1} + \alpha \varpi_1(P_{k+1}) - g_{k} + m_k
\end{multline*}
where the random process $k\to m_k$ is a zero $\{\mathcal{F}_k\}$-martingale comprising all  remaining terms such that $E[m_k|x_k]=0$. By setting
$v_k:=s_k\sqcdot r_{k}$ and $v_{k+1}:=s_{k+1}\sqcdot r_{k+1}$, $k=0,1,\ldots$ we get that 
\begin{multline}\label{eq:dif.one.step}\allowdisplaybreaks
       \alpha^{-k}\left(V^u_{k+1}( x_{k+1}) - V^u_k (x_k)\right)+\|y_k\|^2 =
      x_k^\trp \big(\mathcal{L}^\alpha(P_{k+1})+ C^\trp C -  P_{k} \big) x_k 
      \\
      + \alpha u_k^\trp \Lambda(P_{k+1}) u_k 
       + \langle \alpha A^\trp v_{k+1} + \alpha \mathcal{W}_{x_d}(P_{k+1})s_k - v_{k} ,x_{k} \rangle 
    \\
        +\alpha\langle  B^\trp v_{k+1} + \mathcal{W}_{u}(P_{k+1})s_u +2\Sigma(P_{k+1})x_k ,u_{k} \rangle  
         +\alpha g_{k+1} + \alpha \varpi_1(P_{k+1}) - g_{k} +m_k
\end{multline}
Finally, taking the conditional expectation wrt $x_k$ and using the notation $\mu_k=E[v_k|x_k]$ and $\mu_{k+1}=E[v_{k+1}|x_k]$, we get the expression in the lemma.

\section{Proof of Lemma \ref{lemm:stabilization.G}} \label{app:stabilization.G}

 Let us consider $V^u$ in \eqref{eq:Vu}, and the variation \eqref{eq:dif.one.step.expected} in Lemma \ref{lemm:variation.Vu}. Set $u_k=Gx_k$ and denote $s_x=\mathcal{S}(x_k), s_u=\mathcal{S}(u_k)$ to write,
\begin{multline}\label{eq:dif.one.step.feedback}\allowdisplaybreaks
     E\left[\alpha^{-k}(V^u_{k+1}( x_{k+1}) - V^u_k (x_k))+\|y_k\|^2 |x_k\right]=
     \\
      x_k^\trp \Bigl(\mathcal{L}^\alpha(P_{k+1})+ C^\trp C  
     + \alpha (G^\trp \Lambda(P_{k+1}) G + \Sigma(P_{k+1})^\trp G + G^\trp \Sigma(P_{k+1})\bigl) -  P_{k} \Bigr) x_k 
  \\
       + \bigl\langle \alpha (A+BG)^\trp \mu_{k+1} + \alpha \bigl(\mathcal{W}_{x}(P_{k+1})s_k + G^\trp\mathcal{W}_{u}(P_{k+1})s_u\bigr)  - \mu_{k} ,x_{k} \bigr\rangle 
    \\
      + \alpha g_{k+1} + \alpha \varpi_1(P_{k+1}) - g_{k}
\end{multline}
If there holds,
\begin{multline}\label{eq:quadratic.evol.ctr}
     \mathcal{L}^\alpha(P_{k+1})+ C^\trp C  + 
    \alpha (G^\trp \Lambda(P_{k+1}) G + \Sigma(P_{k+1})^\trp G + G^\trp \Sigma(P_{k+1})\bigl) =  
\\
    \alpha(A+B G)^\trp P_{k+1} (A+B G) + \alpha\mathcal{Z}_x(P_{k+1}) 
+ \alpha G^\trp \mathcal{Z}_u(P_{k+1}) G + (C+DG)^\trp (C+DG)
\\
\mathcal{H}_G(P_{k+1})
=P_{k} 
\end{multline}
this is precisely \eqref{eq:quadratic.G}. Together with \eqref{eq:linear.G} and \eqref{eq:scalar.G}, from \eqref{eq:dif.one.step.feedback} one gets that
\begin{equation*}
    E[V^u_{k}(x_{k}) -V^u_{k+1}(x_{k+1}) | x_k] =  \alpha^{k} E[\| y(k) \|^2|x_k]
\end{equation*}
Since the process $\Theta_{\textrm{ctr}}$ with $u_k=Gx_k$ is Markovian, it follows that
\begin{equation*}
    E[V^u_{0}(x_{0}) - V^u_{\kappa}(x_\kappa)] = E\Bigl[\sum_{k=0}^{\kappa-1} \alpha^k \| y(k) \|^2\Bigr].
 \end{equation*}
and with the choices in the lemma,  $V^u_\kappa(x_\kappa) =0$. This shows the representation part.

To show that $u_k=G_k$ stabilizes $\Theta_{\textrm{ctr}}$, denote $P^{(\kappa)}_{k}, \mu^{(\kappa)}_{k}$ and $g^{(\kappa)}_{k}$, $k=0,\dots,\kappa$, to explicitly indicate the horizon $\kappa$. Note that $0\preceq P_0^{(0)} \preceq P_0^{(1)}\preceq \cdots $ and from the assumption, there exists solution $L\succeq0$ to \eqref{eq:controled.expanded}, 
and hence, $\lim_{\kappa\to\infty}P_0^{(\kappa)}=L$. 

In this situation, from \eqref{eq:quadratic.evol.ctr} we get that  $\alpha^{1/2}(A+BG)$ is d-stable relative to $\alpha(\mathcal{Z}_x+G^\trp\mathcal{Z}_uG)$, which readily implies from Proposition \ref{prop:freiling.stab1}  that  $\Theta_{\textrm{ctr}}$ with $k\to u_k=Gx_k$ is $\alpha$-stable. The $\alpha$-stabilizabilty of $\Theta_{\textrm{ctr}}$ stems straithfowardly.
\end{nolinenumbers}

\end{document}